\documentclass{crm-article}
\usepackage{amsfonts}
\usepackage{amssymb}
\usepackage{amsmath}
\usepackage{mathpazo}
\usepackage{ dsfont } 
\usepackage[mathpazo]{flexisym}
\usepackage{comment}
\usepackage{todonotes}
\usepackage[english]{babel}
\usepackage{secdot}
\usepackage{tikz}
\usepackage{hyperref}
\usepackage{url}
\usepackage{mathtools}
\usepackage{algorithm}
\usepackage{algorithmic}
\usepackage{comment}

\let\argmin\relax

\DeclareMathOperator*{\argmin}{arg\,min}
\renewcommand{\leq}{\leqslant}

\usepackage{nicefrac}

\def \EE {\mathbb E}

\newcommand{\Argmin}{\mathop{\text{Arg}\!\min}}

\def\e{\epsilon}

\begin{document}

%\englishpaper % раскомментировать в том случае, если текст статьи на английском языке

%\norussian % раскомментировать, если нет метаданных на русском

%\affiliationnoref % раскомментировать, если автор один или все авторы из одной организации

%\emailnoref % раскомментировать в том случае, если автор единственный

%\year=2018 % current year by default
\journalVol{10}
\journalNo{1} %выпуска
\setcounter{page}{1}

% раздел журнала
\journalSection{Математические основы и численные методы моделирования}
\journalSectionEn{Mathematical modeling and numerical simulation}

% дата получения
\journalReceived{00.00.0000.}
%\journalReviewed{01.06.2016.}
%принято к публикации
\journalAccepted{00.00.0000.}

\UDC{519.8}
%О некоторых адаптивных ускоренных методах, допускающих альтернированную минимизацию
\title{О связях задач стохастической выпуклой минимизации с задачами минимизации эмпирического риска на шарах в $p$-нормах}
\titleeng{On the relations of stochastic convex optimization problems with empirical risk minimization problems on $p$-norm balls}
\thanks{Исследование выполнено за счет гранта Российского научного фонда (проект № 21-71- 30005).}
\thankseng{This research was funded by Russian Science Foundation (project 21-71- 30005).}

%автор - в формате \author{\firstname{И.\,И.}~\surname{Иванов}}
\author[1,2]{\firstname{Д.\,М.}~\surname{Двинских}}
%автор - в формате \authorfull{Имя Отчество Фамилия}
\authorfull{Дарина Михайловна Двинских}
% автор на англ. - в формате \authoreng{\firstname{I.\,I.}~\surname{Ivanov}}
\authoreng{\firstname{D.\,M.}~\surname{Dvinskikh}}
%автор на англ. - в формате \authorfull{Firstname M. Surname}
\authorfulleng{Darina M. Dvinskikh}
%вписать свою электронную почту
\email{dviny.d@yandex.ru}
%организация - в формате \affiliation{Московский государственный университет,\protect\\ Россия, 141700, г. Москва, ул. Университетская, д. 9}

%автор - в формате \author{\firstname{И.\,И.}~\surname{Иванов}}
\author{\firstname{В.\,В.}~\surname{Пырэу}}
%автор - в формате \authorfull{Имя Отчество Фамилия}
\authorfull{Виталий Вячеславович Пырэу}
% автор на англ. - в формате \authoreng{\firstname{I.\,I.}~\surname{Ivanov}}
\authoreng{\firstname{V.\,V.}~\surname{Pirau}}
%автор на англ. - в формате \authorfull{Firstname M. Surname}
\authorfulleng{Vitali V. Pirau}
%вписать свою электронную почту
\email{pireyvitalik@phystech.edu}
%организация - в формате \affiliation{Московский государственный университет,\protect\\ Россия, 141700, г. Москва, ул. Университетская, д. 9}

%автор - в формате \author{\firstname{И.\,И.}~\surname{Иванов}}
\author[1,2,3]{\firstname{А.\,В.}~\surname{Гасников}}
%автор - в формате \authorfull{Имя Отчество Фамилия}
\authorfull{Александр Владимирович Гасников}
% автор на англ. - в формате \authoreng{\firstname{I.\,I.}~\surname{Ivanov}}
\authoreng{\firstname{A.\,V.}~\surname{Gasnikov}}
%автор на англ. - в формате \authorfull{Firstname M. Surname}
\authorfulleng{Alexander V. Gasnikov}
%вписать свою электронную почту
\email{gasnikov@yandex.ru}
%организация - в формате \affiliation{Московский государственный университет,\protect\\ Россия, 141700, г. Москва, ул. Университетская, д. 9}

% повторите блок для каждого автора;
% если авторов несколько, и автоматическая расстановка сносок от фамилий 
% к организациям приводит к неправильным результатам, укажите правильный 
% вариант в квадраных скобках

\affiliation[1]{141701 Московская обл., Долгопрудный, Институтский пер. 9, 
Московский физико-технический институт (национальный исследовательский университет)}
%организация - в формате \affiliationeng{Moscow State Institute University, 9 University street, Moscow, 141700, Russia}
\affiliationeng[1]{Moscow Institute of Physics and Technology,\protect\\ Dolgoprudny, Russia}
\affiliation[2]{12705, г. Москва, Большой Каретный переулок, д.19 стр. 1, Институт проблем передачи информации РАН им. А.А. Харкевича}
%организация - в формате \affiliationeng{Moscow State Institute University, 9 University street, Moscow, 141700, Russia}
\affiliationeng[2]{Institute for Information Transmission Problems of the Russian Academy of Sciences (Kharkevich Institute)}
\affiliation[3]{385000, Республика Адыгея, г. Майкоп, ул. Первомайская, д. 208,
Кавказский математический центр Адыгейского государственного университета}
%организация - в формате \affiliationeng{Moscow State Institute University, 9 University street, Moscow, 141700, Russia}
\affiliationeng[3]{Caucasus Mathematical Center, Adyghe State University}

%\affiliation[1]{Институт проблем передачи и обработки информации,\protect\\ Россия, Москва}
%\affiliationeng{Institute for Information Transmission Problems RAS,\protect\\ Moscow, Russia}

%\affiliation[1]{Национальный исследовательский университет «Высшая школа экономики»,\protect\\ Россия, Москва}
%\affiliationeng{National Research University Higher School of Economics,\protect\\ Moscow, Russia}

\begin{abstract}
В данной работе рассматриваются задачи выпуклой стохастической оптимизации, возникающие в  анализе данных  (минимизация функции риска), а также в математической статистике (минимизация функции правдоподобия). Такие задачи могут быть решены как онлайн методами, так и оффлайн  (метод Монте-Карло). При оффлайн подходе исходная задача заменяется эмпирической задачей -- задачей минимизации эмпирического риска. В современном машинном обучении ключевым является следующий вопрос: какой размер выборки (количество слагаемых в функционале эмпирического риска) нужно взять, чтобы достаточно точное решение эмпирической задачи  было решением исходной задачи с заданной  точностью.  Базируясь на недавних существенных продвижениях в машинном обучении и оптимизации  для решения выпуклых стохастических задач на евклидовых шарах (или всем пространстве), мы рассматриваем  случай произвольных шаров в $p$-нормах и исследуем как влияет выбор параметра $p$ на оценки необходимого числа слагаемых в функции эмпирического риска. 

В данной работе рассмотрены как выпуклые задачи оптимизации, так и седловые. Для сильно выпуклых задач были обобщены уже  имеющиеся результаты об одинаковых размерах выборки в обоих подходах (онлайн и оффлайн) на произвольные нормы. Более того, было показано, что условие сильной выпуклости может быть ослаблено: полученные результаты справедливы для функций, удовлетворяющих условию квадратичного роста. В случае, когда данное условие не выполняется, предлагается использовать регуляризацию исходной задачи в произвольной норме. В отличие от выпуклых задач, седловые задачи являются намного менее изученными. Для седловых задач размер выборки был получен при условии  $\gamma$-роста седловой функции по разным группам переменных. 
Это условие при $\gamma=1$ есть не что иное, как аналог условие острого минимума в выпуклых задач.  В данной статье было показано, что размер выборки в случае острого минимума (седла) почти не зависит от желаемой точности решения исходной задачи.    
\end{abstract}


\begin{abstracteng}
In this paper, we consider  convex stochastic optimization problems  arising in  machine learning applications  (e.g., risk minimization) and  mathematical statistics (e.g., maximum likelihood estimation).  There are two main approaches  to solve such kinds of problems, namely the Stochastic Approximation approach (online approach)  and the Sample Average Approximation approach, also known as the Monte Carlo approach, (offline approach). In the offline approach, the problem is replaced by its empirical counterpart (the empirical risk minimization problem). The natural question is how to define the problem sample size, i.e., how many realizations should be sampled so that the quite accurate solution of the empirical problem be the solution of the original problem with the desired precision. This issue is one of the main issues in modern machine learning and optimization. In the last decade, a lot 
of significant advances were made in these areas to solve convex stochastic optimization problems on the Euclidean balls (or the whole space). In this work, we are based on these advances and study the case of arbitrary balls in the $\ell_p$-norms. We also explore the question of how  the parameter $p$ affects the estimates of the required number of terms as a function of empirical risk.

In this paper, both convex and saddle point optimization problems are considered. For strongly convex problems, the existing  results on the same sample sizes in both approaches (online and offline) were generalized to arbitrary norms. Moreover, it was shown that the strong convexity condition can be weakened: the obtained results are valid for functions satisfying the quadratic growth condition. In the case when this condition is not met, it is proposed to use the regularization of the original problem in an arbitrary norm. In contradistinction to convex problems, saddle point problems are much less studied. For saddle point problems, the sample size was obtained under the condition of $\gamma$-growth of the objective function.
When $\gamma=1$, this condition is the condition of sharp minimum  in convex problems. In this article, it was shown that the sample size in the case of a sharp minimum is almost independent of the desired accuracy of the solution of the original problem.
\end{abstracteng}
\keywordeng{stochastic optimization, convex optimization, regularization, empirical risk minimization, stochastic approximation, sample average approximation, quadratic growth condition, sharp minimum. }

\maketitle

%Раздел обозначается \paragraph, подраздел - \subparagraph (не \section и \subsection)

\paragraph{Введение} \label{section_1}

Подавляющее число задач математической статистики \cite{spokoiny2015basics,shapiro2021lectures} и машинного обучения \cite{shalev2014understanding,bach2021learning} в конечном итоге сводятся к задачам стохастической оптимизации: минимизации функции риска, представляющей собой математического ожидание функции потерь. Данные задачи можно решать в онлайн режиме \cite{nemirovski2009robust,agarwal2012information} (методами типа стохастического градиентного спуска), когда решение (например, оцениваемый параметр) корректируется по мере поступления новых данных (выборки) и в оффлайн режиме (методом Монте-Карло), когда исходная задача подменяется задачей минимизации функции эмпирического риска \cite{shapiro2005complexity,shalev2009stochastic,shalev2014understanding,bach2021learning} (выборочного среднего функции потерь). Офлайн подход в последние годы стал достаточно популярным в связи с ростом размерностей задач и необходимостью использовать распределенные вычисления \cite{gorbunov2020recent}. Офлайн подход прекрасно позволяет хранить разные части данных (выборки) на разных устройстах. Если онлайн подход для задач выпуклой стохастической оптимизации достаточно хорошо проработан \cite{nemirovski2009robust,agarwal2012information,woodworth2021even}, то в офлайн подходе теоретически обоснованных результатов поменьше \cite{li2021improved}. В частности, если задача стохастической оптимизации рассматривается на неевклидовом шаре, то офлайн подход   не позволяет учитывать такую специфику (за исключением работ \cite{dvinskikh2021stochastic,dvinskikh2021decentralized}, в которых рассматривался один частный случай -- задача на шаре в $1$-норме, без оценок вероятностей больших отклонений), в отличие от онлайн подхода. В настоящей работе устраняется отмеченный недостаток офлайн подхода. 

\paragraph{Основные результаты} \label{section_2}
Рассмотрим задачу стохастической оптимизации
\begin{equation}\label{problem}
  \min_{x \in X} F(x) := \EE_{\xi}{f(x, \xi)}.  
\end{equation}
Как правило, под множеством $X$ будем понимать шар $B_p^d(R)$ радиуса $R$ с центром в точке $0$ в $p$-норме, $p\ge 1$, в пространстве $\mathbb{R}^d$.
\begin{con}[Липшицевость]\label{lipsc} Для всех $x \in X$ и всех $\xi$ выполянется:
    $$|f(y, \xi) - f(x, \xi)| \le M \|y - x\|_p.$$
\end{con}
\begin{con}[Гладкость]\label{smth} Для всех $x \in X$ и всех $\xi$ выполянется:
    $$\|\nabla_x f(y, \xi) - \nabla_x f(x, \xi)\|_q \le L \|y - x\|_p,$$
где $1/p + 1/q = 1$.
\end{con}
% \begin{assumption}[Условие на шум]\label{smth} Существует такая коснтанта $B > 0$, что для всех $x \in X$, всех $\xi$ и всех $k=2,...,N$ выполянется:
%     $$\EE_{\xi}\left[\|\nabla_x f(x_*, \xi) - \nabla_x f(x, \xi)\|_q^k\right] \le \frac{B^{k-2}}{2}k!\EE_{\xi}\left[\|\nabla_x f(x_*, \xi) - \nabla_x f(x, \xi)\|_q^2\right].$$
% \end{assumption}

Задача заключается в определении числа сэмплов (объема выборки) $N$, т.е. независимых одинаково распределенных реализаций случайной величины $\xi$, которое будет достаточно, чтобы некоторый алгоритм (подход) $A$ позволял по $\{\xi^k\}_{k=1}^N$ определить такой $x\left(\{\xi^k\}_{k=1}^N\right)$, что 
\begin{equation}\label{P}
\mathds{P}\left(F\left(x\left(\{\xi^k\}_{k=1}^N\right)\right) - \min_{x\in X}F(x) \le \varepsilon\right)\ge 1 - \sigma.
\end{equation}
Естественно ожидать, что $N$ зависит от $M,L,R,d,\varepsilon,\sigma$. Как будет видно в дальнейшем, существенной зависимости от $L$ в общем случае нет.

Важным местом в приведенном определении является наличие некоторого подхода (алгоритма), обозначенного через $A$, выдающего $x\left(\{\xi^k\}_{k=1}^N\right)$. В действительности оценка параметра $N$ должна также зависеть и от $A$. Принципиально различаются два подхода к тому, как получать $x\left(\{\xi^k\}_{k=1}^N\right)$.

Первый подход -- \textit{онлайн} (в западной литературе часто используется название <<Stochastic Approximation>>). Базируется на процедурах типа (проекции $\pi_X$) \textit{Стохастического градиентного спуска} 
$$x^{k+1} = \pi_X\left\{x^k - h\nabla_x f(x^k,\xi^k)\right\},~~k = 1,...,N$$
и вариациях этого метода \cite{polyak90,polyak1992acceleration,nemirovski2009robust,shapiro2021lectures}. Отметим, что, как правило, в таких процедурах выдается не последняя точка, а среднее по траектории \cite{polyak90}. Большими преимуществами такого подхода является простота получения искомой оценки, возможность адаптивной корректировки оцениваемого вектора параметров $x$ по мере поступления новых данных (выборки). В действительности, именно такие подходы приводят к наилучшим оценкам для параметра $N$ в случае когда $F$ -- выпуклая функция \cite{nemirovski2009robust,agarwal2012information,shapiro2021lectures}.   

Второй подход -- \textit{офлайн}, который также можно называть подходом на основе метода Монте-Карло (в западной литературе часто используется название <<Sample Average Approximation>>) \cite{shapiro2005complexity,shalev2009stochastic,shapiro2021lectures}. 
В основе подхода -- замена функционала задачи \eqref{problem} на выборочное среднее:
\begin{equation}\label{emp_prob}
    \min_{x \in X} \hat{F}(x):= \frac{1}{N}\displaystyle\sum_{k = 1}^{N} f(x, \xi^k).
\end{equation}
Решение (приближенное) задачи \eqref{emp_prob} понимается как $x\left(\{\xi^k\}_{k=1}^N\right)$ при офлайн подходе. Очевидным недостатком подхода является необходимость достаточно точно решать задачу \eqref{emp_prob}. Впрочем, в ряде случаев это может быть и достоинством, если, например, $f(x,\xi)$ обладает дорогим прямым оракулом, выдающим $\nabla_x f(x,\xi)$, но дешевым двойственным, выдающим градиент сопряженной по $x$ функции \cite{dvinskikh2021decentralized,dvinskikh2021stochastic}. Другим недостатком является более скромная теория, которая приводит в целом к худшим оценкам $N$ в выпуклом случае \cite{shapiro2005complexity,shalev2009stochastic,feldman2019high,klochkov2021stability,li2021improved}. Причем (сильную) выпуклость требуется понимать теперь, как (сильную) выпуклость $f(x,\xi)$ по $x$, (сильной) выпуклости только $F$ уже не достаточно  для конечности $N$ \cite{sekhari2021sgd}. Впрочем, как будет видно в дальнейшем, это условие можно заметно ослабить -- в большей степени, чем при онлайн подходе. Отличительным достоинством офлайн подхода является возможность организации распределенных вычислений \cite{gorbunov2020recent} при решении задачи \eqref{emp_prob}, что представляется принципиально важным для многих современных приложений, приходящих, например, из обучения глубоких нейронных сетей \cite{huang2019gpipe}.

Далее в статье постараемся сравнить подробнее оба подхода. Для этого потребуется обобщить некоторые результаты, связанные с офлайн подходом.

% Для большей наглядности изложения будем считать, что хвосты распределения $\xi$ субгауссовские \cite{spokoiny2015basics,shapiro2021lectures}. В действительности, эти предположения можно существенно осл

\subparagraph{Выпуклый случай}\label{cvx}
Для возможности сравнения двух подходов (онлайн и офлайн) предположим, что $f(x,\xi)$ -- удовлетворяет предположению~\ref{lipsc}, а $F(x)=\EE_{\xi}f(x, \xi)$ -- выпуклая функция при $x\in B_p^d(R)$. 

% выпуклая функция по $x$ при $x\in B_p^d(R)$ ($X = B_p^d(R)$) и при всех $\xi$, удовлетворяющая предположению~\ref{lipsc}. Отметим, что для онлайн подхода выпуклость $f(x,\xi)$ можно ослабить до выпуклости $F(x)=\EE_{\xi}f(x, \xi)$ при $x\in B_p^d(R)$.

Из результатов \cite{nemirovski2009robust} следует, что в онлайн подходе 
\begin{itemize}
    \item  при $N \le d$ 
\begin{equation}\label{N_1}
N = O\left(\kappa_p(d)\frac{M^2R^2}{\varepsilon^{\max\{2,p\}}}\ln\left(\frac{1}{\sigma}\right)\right),
\end{equation}
где $\kappa_p(d) = O(1)$, при $p \ge 2$, при $p \in [1,2]$ функция $\kappa_p(d)$ убывает от $O(\ln d)$ при $d = 1$ до $O(1)$ при $p = 2$.
    \item при $N \ge d$ 
\begin{equation}\label{N_2}
N = O\left(d^{1 - 2/\max\{2,p\}}\frac{M^2R^2}{\varepsilon^2}\ln\left(\frac{1}{\sigma}\right)\right).
\end{equation}
\end{itemize}

Причем данные оценки с точностью до логаримических множителей не могут быть улучшены в общем случае, в том числе, даже при дополнительном предположении~\ref{smth} \cite{nemirovksi79,agarwal2012information}.

Из результатов \cite{shapiro2005complexity,shapiro2021lectures} следует, что в офлайн подходе 
\begin{equation}\label{N_erm}
N = O\left(\frac{M^2R^2}{(\varepsilon - \delta)^2}\left(d\ln\left(\frac{MR}{\varepsilon - \delta}\right)+\ln\left(\frac{1}{\sigma}\right)\right)\right),
\end{equation}
где $\delta$ -- точность решения задачи \eqref{emp_prob}. Причем данная оценка с точностью до логарифмических множителей не может быть в общем случае улучшена, в том числе, даже при дополнительном предположении~\ref{smth} \cite{feldman2016generalization}. 

Сопоставляя оценки, которые можно получить при онлайн подходе \eqref{N_1}, \eqref{N_2} с оценкой офлайн подхода, получаем, что за исключением случая $p = \infty$ онлайн подход доминирует офлайн. В частности, при $p=2$ имеем $N_{\text{офлайн}} \simeq d\cdot N_{\text{онлайн}}$.

На самом деле, приведенные выше результаты можно обобщить и на случай, когда $M = M(\xi)$ в предположении~\ref{lipsc} не равномерно ограничена по $\xi$, а ограниченным является лишь второй момент $\EE_{\xi}{\left[M(\xi)^2\right]}$ \cite{gorbunov2021near,shapiro2021lectures}. 

В заключение заметим, что оценка офлайн подхода \eqref{N_erm} может быть получена и без предположения выпуклости функции $F$ \cite{shapiro2021lectures}. То есть выпуклость при офлайн подходе в общем случае ничего не дает. Ситуация существенно меняется в сильно выпуклом случае.

\subparagraph{Сильно выпуклый случай. Условие \textit{квадратичного роста}}\label{scvx}
Отмеченный в 
предыдущем разделе
% разделе~\ref{cvx}
зазор в оценках $N$ в онлайн и офлайн подходе в выпуклом случае исчезает в сильно выпуклом случае \cite{shalev2009stochastic}.

Для простоты сначала предположим, что $f(x,\xi)$ -- $\mu$-сильно выпуклая в $p$-норме функция по $x$ при $x\in X$ ($X$ выпуклое множество) и при всех $\xi$, т.е. для всех $x,y\in X$
\begin{equation}\label{SC}
f(y,\xi)\ge f(x,\xi) + \langle \nabla_x f(x,\xi), y - x \rangle + \frac{\mu}{2}\|y-x\|_p^2.
\end{equation}
Также будем предполагать, что $f(x,\xi)$ удовлетворяет предположению~\ref{lipsc} и является неотрицательной функцией своих аргументов $f(x,\xi)\ge 0$. Отметим, что для онлайн подхода $\mu$-сильную выпуклость в $p$-норме $f(x,\xi)$ можно ослабить до $\mu$-сильной выпуклости в $p$-норме  $F(x)=\EE_{\xi}f(x, \xi)$, а условие неотрицательности $f(x,\xi)$ можно опустить совсем.

Из результатов \cite{juditsky2011first,juditsky2014deterministic,harvey2019tight} следует, что в онлайн подходе 
\begin{equation}\label{N_3}
N = O\left(\kappa_p(d)\frac{M^2}{\mu\varepsilon}\ln\left(\frac{\ln\left( M^2/(\mu\varepsilon)\right)}{\sigma}\right)\right),
\end{equation}
где $\kappa_p(d)$ было определено в
формуле \eqref{N_1}.
% разделе~\ref{cvx}. 
Данная оценка \eqref{N_3} с точностью до логарифмических множителей не может быть улучшена в общем случае, в том числе, даже при дополнительном предположении~\ref{smth} \cite{nemirovksi79}.

Из результатов работ \cite{shalev2009stochastic,feldman2019high,klochkov2021stability,li2021improved}, в которых рассматривался случай $p = 2$, следует, что в офлайн подходе 
\begin{equation}\label{N_erm2}
N = O\left(\frac{M^2}{\mu\varepsilon}\left(\ln\left(\frac{M^2}{\mu\varepsilon}\right) +\ln\ln\left(\frac{1}{\sigma}\right)\right)\ln\left(\frac{1}{\sigma}\right)\right).
\end{equation}
При этом, требуется решить задачу \eqref{emp_prob} с точностью $\delta = O(\mu\varepsilon^2)$. Оценки на $N$ и $\delta$ с точностью до логарифмических множителей не могут быть в общем случае улучшены, в том числе, даже при дополнительном предположении~\ref{smth} \cite{nemirovksi79,shalev2009stochastic}.

В данной работе устанавливается следущий результат.
\begin{teo}\label{MT} Пусть $f(x,\xi)\ge 0$ удовлетворяет условию \eqref{SC} на выпуклом множестве $X$ и удовлетворяет предположению~\ref{lipsc}, где $p\in[1,\infty]$. Пусть задача \eqref{emp_prob}, с $N$ определяемым по формуле \eqref{N_erm2}, решена с точностью по функции $\delta = O(\mu\varepsilon^2)$ с вероятностью $1 - \sigma/2$, т.е. получен такой $x\left(\{\xi^k\}_{k=1}^N\right)$, что 
$$\mathds{P}\left(\hat{F}\left(x\left(\{\xi^k\}_{k=1}^N\right)\right) - \min_{x\in X}\hat{F}(x) \le \delta\right)\ge 1 - \sigma/2.$$ 
Тогда $x\left(\{\xi^k\}_{k=1}^N\right)$ будет $\varepsilon$-решением по функции задачи \eqref{problem} с вероятностью $1 - \sigma$ (см. \eqref{P}), т.е. 
$$\mathds{P}\left(F\left(x\left(\{\xi^k\}_{k=1}^N\right)\right) - \min_{x\in X}F(x) \le \varepsilon\right)\ge 1 - \sigma.$$
\end{teo}
\begin{cor}[условие квадратичного роста]\label{MTcoll}  В условиях теоремы~\ref{MT} можно ослабить
условие сильной выпуклости \eqref{SC} до условия выпуклости $f(x,\xi)$ по $x$ и условия  квадратичного роста функций $\hat{F}(x)$ из \eqref{emp_prob} и $F(x)$ из \eqref{problem}:

для всех $x\in X$ (и всех $\{\xi^k\}_{k=1}^N$, см. \eqref{emp_prob})
\begin{equation}\label{QG}
    \hat{F}(x) - \hat{F}(\hat{x}_*) \ge \frac{\mu}{2}\|x-\hat{x}_*\|_p^2,
\end{equation}
где $\hat{x}_*$ -- проекция $x$ на множество решений задачи \eqref{emp_prob};
\begin{equation}\label{QGF}
    F(x) - F(x_*) \ge \frac{\mu}{2}\|x-x_*\|_p^2,
\end{equation}
где $x_*$ -- проекция $x$ на множество решений задачи \eqref{problem}.
\end{cor}
При $p = 2$ это следствие было установлено в работе \cite{li2021improved}. Также в данной работе приведены другие обобщения приведенной теоремы при $p = 2$, в частности, на случай, когда можно совсем отказаться от условий выпуклости, заменив их намного более слабым условием условием Поляка--Лоясиевича, которму должна удовлетворять функция $F$, а не $\hat{F}$:
 
для всех $x \in X$
\begin{equation}\label{PL}
    F(x) - F(x_*) \le \frac{1}{2\mu}\|\nabla F(x)\|_2^2,
\end{equation}
где $x_*$ -- проекция (в 2-норме) $x$ на множество решений задачи \eqref{problem}.

А именно,  в \cite{li2021improved} показано, что если дополнительно (к условию Поляка--Лоясиевича для $F$) для $f(x, \xi)\ge 0$ выполняются  предположения~\ref{lipsc},~\ref{smth},~\ref{noise}
\begin{con}[Условие на шум]\label{noise} Существует такая константа $B > 0$, что для всех $k=2,...,N$ выполянется:
    $$\EE_{\xi}\left[\|\nabla_x f(x_*, \xi)\|_2^k\right] \le B^{k-2}k!\EE_{\xi}\left[\|\nabla_x f(x_*, \xi) \|_2^2\right],$$
где $x_*$ -- решение задачи \eqref{problem}.
 \end{con}
то при достаточно большом $N$ с вероятностью $1 - \sigma$ справедлива оценка
$$F\left(x\left(\{\xi^k\}_{k=1}^N\right)\right) - F(x_*) = O\left(\frac{(B^2 + \mu^2)\ln^2\left(1/\sigma\right)}{\mu N^2} + \frac{L F(x_*)\ln\left(1/\sigma\right)}{\mu N}\right).$$
В перепараметризованном случае $F(x_*)\simeq 0$ получаем, что $N \sim 1/\sqrt{\mu\varepsilon}$, что сильно лучше оценки \eqref{N_erm2}, но может быть хуже оценки, которую можно получить в перепараметризованном случае для онлайн подхода $N \sim (L/\mu) \ln  \left(\Delta F/\varepsilon\right)$, см., например, \cite{woodworth2021even}.

Интересно было попробовать обобщить и эти результаты на случай $p \in [1,\infty]$. Насколько нам известно, это пока еще не сделано.

\subparagraph{Регуляризация}
Из предыдущих разделов
% разделов~\ref{cvx} и~\ref{scvx}
следует, что в случае выпуклой задачи, выгодно сделать ее сильно выпуклой с помощью \textit{регуляризации} (см., например, \cite{shalev2009stochastic,dvinskikh2021stochastic,dvinskikh2021decentralized}). Причем <<эффект>> от такой регуляризации будет значительно выше, чем это имеет место в обычной оптимизации \cite{nemirovksi79,gasnikov2021}.

Под <<регуляризацией>> понимается замена исходной задачи \eqref{problem} на задачу с $f(x,\xi):= f(x,\xi) + \mu V(x,x^0)$, где $V(x,x^0)$ -- $1$-сильно выпуклая по $x$ на $X$ в $p$-норме ($p\in[1,2]$) функция, такая что (см. обозначения в разделе~\ref{cvx}) $V(x,x^0) \le \kappa_p(d)\|x - x^0\|_p^2  =O\left(\|x - x^0\|_p^2 \ln d \right)$. Можно показать, что такие функции существуют \cite{ben2001lectures}, и уже вполне успешно применялись в рассматриваемом здесь контексте \cite{dvinskikh2021stochastic,dvinskikh2021decentralized} при $p = 1$. В данной работе рассматривается общий случай $p\in[1,2]$. 

Ключевое наблюдение (см., например, замечание 4.1 \cite{gasnikov2021}) заключается в том, что если $x\left(\{\xi^k\}_{k=1}^N\right)$ -- $(\varepsilon/2,\sigma)$-решение регуляризованной задачи в смысле \eqref{P} с $\mu \le \varepsilon/(2V(x_*,x^0))$, где $x_*$ -- такое решение задачи \eqref{problem}, которое наиболее близко к $x^0$ (в смысле минимальности $V(x_*,x^0)$), то $x\left(\{\xi^k\}_{k=1}^N\right)$ будет $(\varepsilon,\sigma)$-решением исходной задачи \eqref{emp_prob} в смысле \eqref{P}.

Выбирая <<на пределе>> $\mu = \varepsilon/\left(2\kappa_p(d)R^2\right) $ получим (с точностью до логаримических множителей) из формул раздела <<Сильно выпуклый случай. ...>>
% ~\ref{scvx} 
формулы раздела <<Выпуклый случай>>,
% ~\ref{cvx},
только без лишнего $d$-множителя в офлайн случае \eqref{N_erm}.

Таким образом, регуляризация решает отмеченную проблему нестыковки оценок онлайн и офлайн подходов в выпуклом случае ($f(x,\xi)$ -- выпуклая функция от $x$). Впервые приблизительно такая конструкция была предложена в данном контексте при $p = 2$ в работе \cite{shalev2009stochastic} (см. также ее изложение, вошедшее в классический учебник по Машинному обучению \cite{shalev2014understanding}), а для $p = 1$ близкая конструкция была описана в работе \cite{dvinskikh2021stochastic}. Описанный выше подход обобщает схему из \cite{dvinskikh2021stochastic} на случай $p \in [1,2]$.

\subparagraph{Промежуточная выпуклость. Острый минимум}\label{sharp}
Условие квадратичного роста можно обобщить. Введем, следуя Шапиро--Немировскому, см., например, \cite{shapiro2005complexity,shapiro2021lectures}, \textit{условие $\gamma$-роста} ($\gamma \ge 1$):

для всех $x\in X_{2\varepsilon} = \left\{x\in X:~~ F(x)\le F(x_*) + 2\varepsilon \right\}$:
\begin{equation}\label{gamma}
    F(x) - F(x_*) \ge \mu_{\gamma}\|x - x_*\|_p^{\gamma},
\end{equation}
где $x_*$ -- проекция (в $p$-норме) $x$ на множество решений задачи \eqref{problem}. 

Ослабим также предположение~\ref{lipsc}. А именно, предположим, что для любых $x,y \in X$ субгауссовская дисперсия $f(y,\xi) - f(x,\xi) - \left( F(y) - F(x)\right)$ ограничена сверху $\lambda^2\|y-x\|_p^2$, т.е.
\begin{equation}\label{var}
\EE_{\xi}\left[\exp\left(t\cdot\left(f(y,\xi) - f(x,\xi) - \left( F(y) - F(x)\right)\right)\right)\right]\le \exp\left(t^2\lambda^2\|y-x\|_p^2/2\right).
\end{equation}
Заметим, что если выполняется предположение~\ref{lipsc}, то $\lambda^2 \le 2 M^2$.

Если $f(x,\xi)$ -- выпуклая по $x$ функция (на самом деле, это условие можно ослабить \cite{shapiro2021lectures}), то при сделанных предположениях\footnote{Если $M$ в предположении~\ref{lipsc} зависит от $\xi$, то под $M$ в формуле \eqref{N_erm3} следует понимать $\EE_{\xi} M(\xi)$ \cite{shapiro2005complexity,shapiro2021lectures}. Параметр $R_{\varepsilon}$ в этой формуле отвечает диаметру множества $X_{2\varepsilon}$ в $p$-норме. В частности, при $\gamma = 1$ параметр $R_{\varepsilon}\le 4\varepsilon/\mu_1$. Таким образом, в случае <<острого минимума>> ($\gamma = 1$) $N$ не зависит от $\varepsilon$ \cite{shapiro2021lectures}.}
\begin{equation}\label{N_erm3}
N = O\left(\frac{\lambda^2}{\mu_{\gamma}^{2/\gamma}\varepsilon^{2(\gamma - 1)/\gamma}}\left(d\ln\left(\frac{MR_{\varepsilon}}{\varepsilon}\right)+\ln\left(\frac{1}{\sigma}\right)\right)\right),
\end{equation}
где $\delta = \varepsilon/2$ -- точность решения задачи \eqref{emp_prob}. Причем данная оценка \eqref{N_erm3} с точностью до логарифмических множителей не может быть в общем случае улучшена \cite{shapiro2005complexity,shapiro2021lectures}. В цитированных работах оценка \eqref{N_erm3} была доказана, насколько удалось понять обозначения, для случая $p=2$. Однако, в \cite{shapiro2005complexity,shapiro2021lectures} общий случай $p\in[1,\infty)$ получается дословным повторением всех рассуждений, что также было нам подтверждено в ходе личной беседы одним из авторов \cite{shapiro2005complexity,shapiro2021lectures} Александром Шапиро. 

Формула \eqref{N_erm3} особенно интересна в случае <<острого минимума>> $\gamma = 1$. Она не зависит от $\varepsilon$.

То что формула \eqref{N_erm3} не может быть улучшена хорошо поясняет пример из книги \cite{shapiro2021lectures}, в котором $p = 2$: $f(x,\xi) = \|x\|_2^{\gamma} - \gamma\sigma\langle\xi,x\rangle$, $\xi \in \mathcal{N}(0,I_d)$ -- стандартное нормальное распределение (с нулевым математическим ожиданием и единичной корреляционной матрицей), $X = B_2^d(1)$. В этом случае $N$ не может быть меньше чем $d\sigma^2/\varepsilon^{2(\gamma-1)/\gamma}$. Однако мы привели здесь этот пример, чтобы показать, что предположение~\ref{lipsc} и условие~\eqref{var} могут довольно сильно отличаться. А именно, для этого примера предположение~\ref{lipsc} выполняется лишь в <<среднем>> с $\EE_{\xi}{\left[M(\xi)^2\right]} = \gamma^2\sigma^2d$, притом, что условие~\eqref{var} выполняется с $\lambda = \gamma^2\sigma^2$.

В связи со всем выше написанным в этом разделе и написанным ранее в разделе
<<Выпуклый случай>>
% ~\ref{cvx} 
может показаться, что оценка \eqref{N_erm3} при $\gamma \to \infty$ (вырожденный случай) противоречит нижней оценке \eqref{N_erm} \cite{feldman2016generalization}. Ведь оценка $N$ сверху \eqref{N_erm3} получается лучше в плане возможности использования параметра $\lambda$ вместо $M$. Но при этом, предположение~\ref{lipsc} является более узким, чем условие~\eqref{var}. На самом деле, никакого противоречия нет. Обе оценки точные в своих классах функций $f(x, \xi)$. Возникающий здесь парадокс с описанным примером связан с тем, что $\gamma \to \infty$ влечет за собой то, что $\lambda \to \infty$ и $M \to \infty$. Поэтому данный пример в пределе $\gamma \to \infty$ не отражает точное поведение оценки \eqref{N_erm3}. 

Результат, аналогичный \eqref{N_erm3} (с заменой $d\lambda^2$ на $M^2$ и $R_{\e}$ на $R$ при $\gamma = 1$) при условии~\ref{gamma} может быть получен и для онлайн методов типа рестартованного стохастического градиентного спуска  при $\gamma \ge 2$ \cite{juditsky2014deterministic}. Результаты работы \cite{juditsky2014deterministic} переносится и на случай $\gamma \in [1,2]$. Случай $\gamma = 1$ был исследован в работе \cite{juditsky1993}.

\subparagraph{Седловые задачи. Промежуточная выпукло-вогнутость. Острый минимум}\label{sharp_interm}
К сожалению, такой богатой теории, которая уже создана для задач (выпуклой) оптимизации, для седловых задач нам не известно. Из всех приведенных выше результатов на данный момент удалось перенести только результат \eqref{N_erm3}. Рассуждения практически дословно повторяют выкладки из работ \cite{shapiro2005complexity,shapiro2021lectures}. Далее излагается соответствующая теория.

Рассматривается стохастическая седловая задача
\begin{equation}\label{saddleproblem}
  \min_{x \in X} \max_{y \in Y} F(x,y) := \EE_{\xi}{f(x,y, \xi)}.  
\end{equation}
Множества $X\subset \mathbb{R}^{d_x}$ и $Y\subset \mathbb{R}^{d_y}$ предполагаются выпуклыми компактами, функция $f(x,y,\xi)$ -- выпуклая по $x$ при $x \in X$ и вогнутая по $y$ при $y \in Y$ для всех $\xi$. Также будем считать, что по каждой группе переменных $x$ и $y$ на $X$ и на $Y$ функция $f(x,y,\xi)$ удовлетворяет: 
\begin{itemize}
    \item предположению \ref{lipsc} c параметрами, соответственно, $M_x(\xi)$ в $p_x$-норме ($p_x \in [1,2]$) и $M_y(\xi)$ в $p_y$-норме ($p_y \in [1,2]$), причем, $\EE M_x(\xi) = M_x < \infty$ и $\EE M_y(\xi) = M_y < \infty$;
    \item условию \eqref{var} с параметрами, соответственно, $\lambda_x$, $p_x$ и $\lambda_y$, $p_y$;
    \item условию $\gamma_x$-роста в $p_x$-норме с константой $\mu_{\gamma,x}$ и $\gamma_y$-роста в $p_y$-норме с константой $\mu_{\gamma,y}$:
\begin{equation}\label{gamma_x}
    F(x,y) - F(x_*(y),y) \ge \mu_{\gamma,x}\|x - x_*(y)\|_{p_x}^{\gamma_x},
\end{equation}
где $x_*(y)$ -- проекция (в $p_x$-норме) $x$ на множество решений задачи  $\min_{x\in X}  F(x,y)$ и 
\begin{equation}\label{gamma_y}
   F(x,y_*(x)) - F(x,y) \ge \mu_{\gamma,y}\|y - y_*(x)\|_{p_y}^{\gamma_y},
\end{equation}
где $y_*(x)$ -- проекция (в $p_y$-норме) $y$ на множество решений задачи  $\max_{y\in Y} F(x,y)$.
\end{itemize}

Введем эмпирическую функцию
$$\hat{F}(x,y) = \frac{1}{N}\sum_{k=1}^N f(x,y,\xi^k).$$
Путь удалось найти такие $(\hat{x},\hat{y})$, что 
$$ \hat{F}(\hat{x},\hat{y}) - \hat{F}(\hat{x}_*(\hat{y}),\hat{y})\le \varepsilon/4,$$
$$ \hat{F}(\hat{x},\hat{y}_*(\hat{x}))- \hat{F}(\hat{x},\hat{y})  \le \varepsilon/4,$$
где $\hat{x}_*(y)$, $\hat{y}_*(x)$ определяются по $\hat{F}$ аналогично тому, как $x_*(y)$, $y_*(x)$ определялись по $F$.

Тогда если
$$
N = O\left(N_x + N_y \right),
$$
$$
N_x = O\left(\frac{\lambda_x^2}{\mu_{\gamma,x}^{2/\gamma_x}\varepsilon^{2(\gamma_x - 1)/\gamma_x}}\left(d_x\ln\left(\frac{M_xR_x}{\varepsilon}\right)+\ln\left(\frac{1}{\sigma}\right)\right)\right),
$$
$$
N_y = O\left(\frac{\lambda_y^2}{\mu_{\gamma,y}^{2/\gamma_y}\varepsilon^{2(\gamma_y - 1)/\gamma_y}}\left(d_y\ln\left(\frac{M_yR_y}{\varepsilon}\right)+\ln\left(\frac{1}{\sigma}\right)\right)\right),
$$    
то с вероятностью $1 - \sigma$
$$ F(\hat{x},\hat{y}) - F(x_*(\hat{y}),\hat{y})\le \varepsilon/2,$$
$$ F(\hat{x},y_*(\hat{x})) - F(\hat{x},\hat{y}) \le \varepsilon/2,$$
где $R_x$ -- диаметр $X$ в $p_x$-норме, $R_y$ -- диаметр $Y$ в $p_y$-норме (можно уточнить эти оценки и использовать диаметры соответствующих множеств Лебега, подобно тому, как это делалось выше, см. сноску 1). Следовательно,
$$0\le \max_{y\in Y}F(\hat{x},y) - \min_{x\in X} F(x,\hat{y}) = F(\hat{x},y_*(\hat{x})) - F(x_*(\hat{y}),\hat{y}) \le \varepsilon.$$

Приведенные выше оценки могут быть получены и в онлайн подходе (с заменой $d\lambda^2$ на $M^2$). Немного в более общем контексте это недавно было показано в работе \cite{dvinskikh2022}.

Авторы выражают благодарность Александру Шапиро и Анатолию Юдицкому за ценные советы.

Статья приурочена к 60-и летию Анатолия Борисовича Юдицкого, внесшего значительный вклад в развитие методов стохастического градиентного спуска.


\begin{thebibliography}{99}
\bibitem[Гасников, 2021]{gasnikov2021}Гасников, А. Современные численные методы оптимизации. Метод универсального градиентного спуска. (М.: МЦНМО,2021)\\
Gasnikov A.V. Sovremennye chislennye metody optimizatsii. Metod universal'nogo gradientnogo spuska. [Universal gradient method] MCCME, 2021. (in Russian)
\bibitem[Немировский and Юдин, 1979]{nemirovksi79}Немировский, А.C. \& Юдин, Д.Б. Сложность задач и эффективность методов оптимизации. (Наука,1979)\\
Nemirovsky A.S., Yudin D.B. Slozhnost' zadach i effektivnost' metodov optimizatsii. [Problem Complexity and Optimization Method Efficiency]. M.: Nauka. – 1979.(in Russian)
\bibitem[Поляк, 1990]{polyak90}Поляк, Б.Т. Новый метод типа стохастической аппроксимации. {\em Автоматика и Телемеханика}., 98-107 (1990)\\
Polyak, B.T.  Novyi metod tipa stokhasticheskoi approksimatsii. [A new method of stochastic approximation type”, Autom. Remote Control], 51:7 (1990), 937–946.
\bibitem[Agarwal et al., 2012]{agarwal2012information}Agarwal, A., Bartlett, P., Ravikumar, P. \& Wainwright, M. Information-theoretic lower bounds on the oracle complexity of stochastic convex optimization.// Advances in Neural Information Processing Systems., 2009., V. \textbf{22}.
%{\em IEEE Transactions On Information Theory}. \textbf{58}, 3235-3249 (2012)
\bibitem[Bach, 2021]{bach2021learning}Bach, F. Learning Theory from First Principles Draft.  (2021)
\bibitem[Ben-Tal and Nemirovski, 2022]{ben2001lectures}Ben-Tal, A. \& Nemirovski, A. Lectures on modern convex optimization: analysis, algorithms, and engineering applications. (SIAM, 2022),\\ \url{https://www2.isye.gatech.edu/~nemirovs/LMCOLN2022WithSol.pdf}
\bibitem[Boucheron et al., 2013]{boucheron2013concentration}Boucheron, S., Lugosi, G. \& Massart, P. Concentration inequalities: A nonasymptotic theory of independence. (Oxford university press,2013)
\bibitem[Dvinskikh, 2021a]{dvinskikh2021decentralized}Dvinskikh, D. Decentralized Algorithms for Wasserstein Barycenters. {\em ArXiv Preprint ArXiv:2105.01587}. Дис. – Humboldt Universitaet zu Berlin (Germany). (2021)
\bibitem[Dvinskikh, 2022]{dvinskikh2022}
Dvinskikh D. et al. Gradient-Free Optimization for Non-Smooth Minimax Problems with Maximum Value of Adversarial Noise. {\em ArXiv Preprint ArXiv:2202.06114.} (2022)
\bibitem[Dvinskikh, 2021b]{dvinskikh2021stochastic}Dvinskikh, D. Stochastic approximation versus sample average approximation for Wasserstein barycenters. {\em Optimization Methods And Software}. pp. 1-33 (2021)
\bibitem[Feldman, 2016]{feldman2016generalization}Feldman, V. Generalization of erm in stochastic convex optimization: The dimension strikes back. {\em Advances In Neural Information Processing Systems}. V. \textbf{29} pp. 3576-3584, (2016)
\bibitem[Feldman and Vondrak, 2019]{feldman2019high}Feldman, V. \& Vondrak, J. High probability generalization bounds for uniformly stable algorithms with nearly optimal rate. {\em Conference On Learning Theory},  PMLR, 2019., pp. 1270-1279 (2019)
\bibitem[Gorbunov et al., 2021]{gorbunov2021near}Gorbunov, E., Danilova, M., Shibaev, I., Dvurechensky, P. \& Gasnikov, A. Near-Optimal High Probability Complexity Bounds for Non-Smooth Stochastic Optimization with Heavy-Tailed Noise. {\em ArXiv Preprint ArXiv:2106.05958}. (2021)
\bibitem[Gorbunov et al., 2020]{gorbunov2020recent}Gorbunov, E., Rogozin, A., Beznosikov, A., Dvinskikh, D. \& Gasnikov, A. Recent theoretical advances in decentralized distributed convex optimization. {\em ArXiv Preprint ArXiv:2011.13259}. (2020)
\bibitem[Harvey et al., 2019]{harvey2019tight}Harvey, N., Liaw, C., Plan, Y. \& Randhawa, S. Tight analyses for non-smooth stochastic gradient descent. {\em Conference On Learning Theory}, PMLR, 2019., pp. 1579-1613 (2019)
\bibitem[Huang et al., 2019]{huang2019gpipe}Huang, Y., Cheng, Y., Bapna, A., Firat, O., Chen, D., Chen, M., Lee, H., Ngiam, J., Le, Q., Wu, Y. \& Others Gpipe: Efficient training of giant neural networks using pipeline parallelism. {\em Advances In Neural Information Processing Systems}, V. \textbf{32} pp. 103-112 (2019)
\bibitem[Juditsky, 1993]{juditsky1993}
Juditsky A. A stochastic estimation algorithm with observation averaging. {\em IEEE transactions on automatic control.} 1993. V. 38. no. 5. P. 794--798.
\bibitem[Juditsky et al., 2011]{juditsky2011first}Juditsky, A., Nemirovski, A.  First order methods for nonsmooth convex large-scale optimization, i: general purpose methods. {\em Optimization For Machine Learning}. V. \textbf{30} 30., no. 9, 121-148 (2011)
\bibitem[Juditsky and Nesterov, 2014]{juditsky2014deterministic}Juditsky, A. \& Nesterov, Y. Deterministic and stochastic primal-dual subgradient algorithms for uniformly convex minimization. {\em Stochastic Systems}. V. \textbf{4}. – no. 1., 44-80 (2014)
\bibitem[Klochkov and Zhivotovskiy, 2021]{klochkov2021stability}Klochkov, Y. \& Zhivotovskiy, N. Stability and Deviation Optimal Risk Bounds with Convergence Rate  $O (1/n)$. {\em ArXiv Preprint ArXiv:2103.12024} // Advances in Neural Information Processing Systems., 2021., V. \textbf{34} (2021)
\bibitem[Li and Liu, 2021]{li2021improved}Li, S. \& Liu, Y. Improved Learning Rates for Stochastic Optimization: Two Theoretical Viewpoints. {\em ArXiv Preprint ArXiv:2107.08686}. (2021)
\bibitem[Nemirovski et al., 2009]{nemirovski2009robust}Nemirovski, A., Juditsky, A., Lan, G. \& Shapiro, A. Robust stochastic approximation approach to stochastic programming. {\em SIAM Journal On Optimization}. V. \textbf{19}, no. 4, 1574-1609 (2009)
\bibitem[Polyak and Juditsky, 1992]{polyak1992acceleration}Polyak, B. \& Juditsky, A. Acceleration of stochastic approximation by averaging. {\em SIAM Journal On Control And Optimization}. V. \textbf{30}, no.4, 838-855 (1992)
\bibitem[Robbins and Monro, 1951]{robbins1951stochastic}Robbins, H. \&  Monro, S. A stochastic approximation method. {\em The annals of mathematical statistics}. V. \textbf{2} pp. 400--407 (1951)
\bibitem[Sekhari et al., 2021]{sekhari2021sgd}Sekhari, A., Sridharan, K. \& Kale, S. SGD: The Role of Implicit Regularization, Batch-size and Multiple-epochs. {\em Advances In Neural Information Processing Systems}. V. \textbf{34} (2021)
\bibitem[Shalev-Shwartz and Ben-David, 2014]{shalev2014understanding}Shalev-Shwartz, S. \& Ben-David, S. Understanding machine learning: From theory to algorithms. (Cambridge university press, 2014)
\bibitem[Shalev-Shwartz et al., 2009]{shalev2009stochastic}Shalev-Shwartz, S., Shamir, O., Srebro, N. \& Sridharan, K. Stochastic Convex Optimization.. {\em COLT}. V.,\textbf{2}, no. 4., (2009)
\bibitem[Shapiro et al., 2021]{shapiro2021lectures}Shapiro, A., Dentcheva, D. \& Ruszczynski, A. Lectures on stochastic programming: modeling and theory. (SIAM, 2021)
\bibitem[Shapiro and Nemirovski, 2005]{shapiro2005complexity}Shapiro, A. \& Nemirovski, A. On complexity of stochastic programming problems. {\em Continuous Optimization}, Springer, Boston, pp. 111-146 (2005)
\bibitem[Spokoiny and Dickhaus, 2015]{spokoiny2015basics}Spokoiny, V. \& Dickhaus, T. Basics of modern mathematical statistics. (Heidelberg, Springer, 2015)
\bibitem[Woodworth and Srebro, 2021]{woodworth2021even}Woodworth, B. \& Srebro, N. An Even More Optimal Stochastic Optimization Algorithm: Minibatching and Interpolation Learning // Advances in Neural Information Processing Systems. – 2021. – V. 34. {\em ArXiv Preprint ArXiv:2106.02720}. (2021)
\end{thebibliography}
\end{document}